\documentclass[12pt,francais,nothms]{aclart}
\usepackage{mathptm}
\usepackage[all]{xy}

\theoremstyle{remark}

\usepackage[all]{xy}

\RequirePackage[T1]{fontenc}
\RequirePackage{mathrsfs,amsmath,amssymb}
\let\mathcal\mathscr
 
\RequirePackage{url}

\RequirePackage{xspace}
\RequirePackage[english,frenchb]{babel}

\theoremstyle{plain}
\newtheorem{lemm}{Lemme}
\newtheorem*{theo}{Th\'eor\`eme}
\theoremstyle{remark}
\newtheorem*{rema}{Remarque}

\let\ra\rightarrow
\def\P{\mathbf P}
\def\C{\mathbf C}
\def\Spec{\operatorname{Spec}}
\def\cf{\emph{cf.}~\ignorespaces}
\def\resp{\emph{resp.}~\ignorespaces}
\def\im{\operatorname{im}}

\begin{document}

\selectlanguage{french}
\title[Groupe fondamental des vari\'et\'es rationnellement connexes]
  {{\normalsize\bfseries G\'eom\'etrie alg\'ebrique/\itshape Algebraic geometry \\}
   \`A propos du groupe fondamental \\
   des vari\'et\'es rationnellement connexes}
\alttitle{On the fundamental group of rationnally connected varieties}
\author{Antoine Chambert-Loir}
\address {Centre de math\'ematiques (CMAT) \\ \'Ecole polytechnique \\ 
91128 Palaiseau Cedex}
\email{chambert@math.polytechnique.fr}

\date{5 mars 2003}

\begin{altabstract}
I show that the fundamental group of a normal variety which is rationally
chain connected is finite. The proof holds in non-zero characteristic.

{\it To cite this article: A. Chambert-Loir,
C. R.  Acad. Sci. Paris, Ser. I ??? (2003).}
\end{altabstract}

\begin{abstract}
Je d\'emontre que le groupe fondamental d'une vari\'et\'e normale rationnellement
connexe par cha\^{\i}nes est fini. La d\'emonstration est valable en
caract\'eristique diff\'erente de z\'ero.

{\it Pour citer cet article~: A. Chambert-Loir, C. R.
Acad. Sci. Paris, Ser. I ??? (2003).}
\end{abstract}

\maketitle

\selectlanguage{french}

Cette note a pour but de donner la d\'emonstration
du th\'eor\`eme suivant :

\begin{theo}\label{theo.fini}
Soit $k$ un corps alg\'ebriquement clos et soit $X$ un $k$-sch\'ema
normal qui est rationnellement connexe par cha\^{\i}nes.

\emph{a}) Le groupe fondamental de~$X$ est fini ;

\emph{b)} Si $k=\C$, il en est de m\^eme de son groupe fondamental
topologique.

\emph{c)} Si la caract\'eristique de~$k$ est un nombre premier~$p$,
le groupe fondamental de~$X$ est d'ordre premier \`a~$p$.
\end{theo}

Rappelons qu'une $k$-vari\'et\'e est dite rationnellement connexe par
cha\^{\i}nes si pour tout corps alg\'ebriquement clos $\Omega$
contenant~$k$, il passe par deux points g\'en\'eraux de~$X(\Omega)$ une courbe
connexe propre dont la normalis\'ee est une r\'eunion disjointe de droites
projectives --- autrement dit, deux points g\'en\'eraux sont joints
par une cha\^{\i}ne de courbes rationnelles propres.

Si $k$ est de caract\'eristique nulle et $X$ lisse, 
on sait plus g\'en\'eralement
que $X$ est alg\'ebriquement simplement connexe, et m\^eme topologiquement
simplement connexe si $k=\C$ 
(Campana~\cite{campana1994}, Koll\'ar, Miyaoka, Mori~\cite{kollar-m-m1992}).
La d\'emonstration de Koll\'ar, Miyaoka et Mori repose sur le fait que des
cha\^{\i}nes de courbes rationnelles peuvent \^etre contract\'ees en
une seule (lissification de courbes rationnelles),
et cela n\'ecessite que $X$ soit propre, lisse 
et d\'efinie sur un corps de caract\'eristique z\'ero.
Celle de Campana ne lissifie pas les courbes rationnelles
mais utilise la finitude du groupe fondamental
topologique et la th\'eorie
de Hodge qui, sous ces hypoth\`eses,  entra\^{\i}ne
que $H^i(X,\mathcal O_X)=0$ pour $i>0$.

Ce dernier point n'est pas n\'ecessairement v\'erifi\'e si la caract\'eristique
du corps de base est diff\'erente de z\'ero, comme le montre l'exemple
des surfaces construites par Shioda~\cite{shioda1974} qui sont \`a la fois
de type g\'en\'eral et unirationnelles ; elles ne sont d'ailleurs pas simplement
connexes. En revanche, il r\'esulte
du th\'eor\`eme de de Jong et Starr~\cite{dejong-s2003} et d'arguments
de Koll\'ar (voir~\cite{debarre2002})
qu'une vari\'et\'e \emph{s\'eparablement rationnellement connexe}
est simplement connexe.
Je ne sais pas si la m\'ethode que suit Koll\'ar
dans~\cite{kollar1995} pour d\'emontrer le th\'eor\`eme~\ref{theo.fini}
en caract\'eristique z\'ero permet de traiter le cas g\'en\'eral.

Rappelons aussi que d'apr\`es un th\'eor\`eme fondamental 
de Campana~\cite{campana1992}
et Kolll\'ar, Miyaoka, Mori~\cite{kollar-m-m1992b},
les vari\'et\'es de Fano
(propres, lisses, connexe \`a fibr\'e anticanonique ample)
sont rationnellement connexes par cha\^{\i}nes. Leur d\'emonstration
est valable en caract\'eristique quelconque (\cf\cite{debarre2001}).
Le groupe fondamental d'une vari\'et\'e de Fano est donc toujours fini.

\section*{D\'emonstration}

\begin{lemm}\label{lem.2}
Soit $\Omega$ une extension alg\'ebriquement close de~$k$.
Soit $F_\Omega\colon \P^1_\Omega\ra X_\Omega$
une courbe rationnelle de $X_\Omega$. Soit $x_0=F_\Omega(0)$
et $x_\infty=F_\Omega(\infty)$ ; ce sont des points de $X_\Omega$.
Notons $V_0$ et $V_\infty$ l'adh\'erence de leurs images dans~$X$.

Alors, il existe un $k$-sch\'ema int\`egre $T$, un morphisme
$F\colon \P^1_T\ra X$ tel que les morphismes
$F_0$ et $F_\infty$ d\'efinis par $F_0(t)=F(0,t)$ et $F_\infty=F(\infty,t)$
soient dominants sur $V_0$ et $V_\infty$.
\end{lemm}
Il existe une $k$-alg\`ebre de type fini~$A$, contenue dans~$\Omega$
et un morphisme $F_A\colon \P^1_A\ra X_A$ tel que
$F_\Omega=F_A\otimes_A \Omega$.
Soit $F$ le corps des fractions de~$A$.
On a alors le diagramme de sch\'emas
\[ \xymatrix{
     \Spec A \ar @<2pt> [r]^-{0} \ar @/^3ex/[rr]^{F_0}
             \ar @<-2pt> [r]_-\infty \ar @/_3ex/[rr]_{F_\infty}
      & \P^1_A \ar[r]^-{F_A} & X_A \ar[r]^{p} & X }
\]
Si $F$ est le corps des fractions de $A$, cela signifie
que les points $x_0$ et $x_\infty$ appartiennent \`a $X(F)$,
et qu'ils s'\'etendent en des morphismes $f_0$, $f_\infty\colon \Spec A\ra X$.
L'image par~$f_0$ du 
le point g\'en\'erique de $\Spec A$ est le point g\'en\'erique
de $V_0$. 
Par suite, l'image de $f_0$ contient un ouvert dense de $V_0$,
et de m\^eme pour $f_\infty$.
(Comme $f_0$ est continue, $f_0^{-1}(V_0)$ est une partie
ferm\'ee de $\Spec A$ qui contient son point g\'en\'erique (d'image $x_0$),
donc est \'egale \`a $\Spec A$.
R\'eciproquement, $f_0(\Spec A)$ est une partie constructible de~$X$ 
qui contient $x_0$, donc elle contient un ouvert dense de $V_0$.)

\bigskip

Pour d\'emontrer le th\'eor\`eme, nous pouvons supposer que le corps~$k$
n'est pas d\'enombrable. Soit $d$ un entier, $0\leq d\leq \dim X$.
L'ensemble des sous-sch\'emas int\`egres de~$X$
qui sont de dimension~$d$ n'est pas d\'enombrable, si bien
que toute famille d\'enombrable de sous-vari\'et\'es de~$X$
\'evite au moins un point de~$X$ dont l'adh\'erence 
est de dimension~$d$.
Si $\Omega$ d\'esigne la cl\^oture alg\'ebrique du corps des fonctions de $X$,
il existe alors dans $X(\Omega)$ des points g\'en\'eraux dont
l'adh\'erence dans~$X$ est de toute dimension, en particulier
de dimension~$0$ (un point de $X(k)$) ou de dimension~$\dim X$
(un point {\og g\'en\'erique\fg}).

La d\'efinition que $X$ est rationnellement connexe par cha\^{\i}nes
implique qu'il existe une cha\^{\i}ne de courbes rationnelles
qui joint un point g\'en\'eral $x_0\in X(k)$
\`a un point g\'en\'erique $x_m\in X(\Omega)$. 
D'apr\`es le lemme~\ref{lem.2},
il existe une suite de sous-vari\'et\'es int\`egres
$V_0,\dots,V_m\subset X$, $V_0$ \'etant un point,
et $V_m=X$, et pour tout entier $i\in\{0;\dots;m-1\}$
une famille de courbes rationnelles int\`egres
\[ \mathbf P^1 \times T \xrightarrow F X \]
telle que les morphismes
$F_0\colon T\ra X$ et $F_\infty\colon T\ra X$,
d\'efinis par $F_0(t)=F(0,t)$ et $F_\infty(t)=F(\infty,t)$
soient d'image denses dans $V_i$ et $V_{i+1}$ respectivement.

Quitte \`a remplacer $T$ par un ouvert, on peut supposer qu'il est
normal, ainsi que $V_i$ et $V_{i+1}$.

\begin{lemm}
\label{lem.3}
Soit $f\colon X\ra Y$ un morphisme dominant
entre sch\'emas int\`egres normaux. Soit $x$ un point
g\'eom\'etrique de~$f$, $y=f(x)$. L'image de $\pi_1(X,x)$
dans $\pi_1(Y,y)$ est d'indice fini.
\end{lemm}

D'apr\`es ce lemme, l'image de $\pi_1(T)$ 
dans $\pi_1(V_i)$ et dans $\pi_1(V_{i+1})$
est d'indice fini. Le diagramme
\[ \xymatrix{
   \pi_1(T) \ar[r]  \ar[d]_{f_0}^{\sim} & \pi_1(V_i) \ar[rd] \\
  \pi_1(\P^1\times T) \ar[rr] && \pi_1(X) \\
   \pi_1(T) \ar[r] \ar[u]_{\sim}^{f_\infty} & \pi_1(V_{i+1}) \ar[ru]}
\]
montre alors que les images de  $\pi_1(V_i)$ et $\pi_1(V_{i+1})$
dans~$\pi_1(X)$ sont commensurables, le sous-groupe
image de $\pi_1(\P^1_T)$ \'etant d'indice fini dans chacun d'eux.

Comme $V_0=\{x_0\}$ et $V_m=X$, le groupe fondamental de~$X$ est
fini, ce qui termine la d\'emonstration des assertions \emph{a})
et \emph{b}) du th\'eor\`eme.

On d\'emontre alors l'assertion \emph{c)} par des techniques
de cohomologie $p$-adique qui combinent
le th\'eor\`eme d'Esnault~\cite{esnault2003}
et des arguments d'Ekedahl~\cite{ekedahl1983}.
Je renvoie \`a~\cite{acl2003b} pour les d\'etails.

\bigskip
\noindent{\itshape
D\'emonstration du lemme~\ref{lem.3}. ---}
P.~Deligne d\'emontre ce lemme \`a la fin de \emph{Th\'eorie de Hodge II},
dans le cas de sch\'emas lisses sur~$\C$
(\cite{deligne1972}, lemme 4.4.17).
Sa d\'emonstration reste valable quasiment mot pour mot 
dans le cas qui nous int\'eresse.

Soit en effet $t$ un point ferm\'e de la fibre g\'en\'erique de~$f$,
$T$ son adh\'erence de Zariski de~$f$. C'est un sous-sch\'ema ferm\'e
de~$X$ et $f$ induit un morphisme g\'en\'eriquement fini $f\colon T\ra Y$.
Il existe un ouvert dense $U\subset Y$ tel que $V=f^{-1}(U)\cap T$
soit un rev\^etement \'etale de~$U$. 
Consid\'erons alors le diagramme de groupes fondamentaux
\[ \xymatrix{
        \pi_1(V) \ar[r] \ar[d]^{\alpha} & \pi_1(X) \ar[d]^{f_*} \\
       \pi_1(U) \ar[r]^{\beta} & \pi_1(Y).  }
\]
L'image de $\pi_1(V)$ par l'application $\alpha$ 
est d'indice fini dans~$\pi_1(U)$ car $V\ra U$ est un rev\^etement \'etale.

D'autre part, $U$ est un ouvert dense dans~$Y$ qui est un sch\'ema
normal ; par suite, l'homomorphisme $\pi_1(U)\ra \pi_1(Y)$ est
surjectif. Cela signifie en effet que la restriction \`a~$U$
d'un rev\^etement \'etale connexe, \resp d'un rev\^etement topologique
connexe de~$Y$ est connexe. 
Consid\'erons un tel rev\^etement $\pi\colon\tilde Y\ra Y$.
Si $Z=X\setminus U$, on a
$\pi^{-1}(U)=\tilde Y\setminus \pi^{-1}(Z)$ 
et $\pi^{-1}(Z)$ est un ferm\'e alg\'ebrique strict
(\resp un espace analytique {\og nulle part dense\fg}
dans le cas topologique) de~$\tilde Z$.
De plus, $\tilde Y$ est un sch\'ema normal  
(\resp un espace analytique normal dans le cas topologique).
donc est irr\'eductible. 
Dans le cas alg\'ebrique, $\pi^{-1}(U)$ est alors irr\'eductible,
donc connexe. Dans le cas topologique,
il r\'esulte du th\'eor\`eme p.~145 de~\cite{grauert-r1984}
que $\pi^{-1}(U)$ est connexe.

L'image de $\pi_1(V)$ dans~$\pi_1(Y)$ est ainsi
d'indice fini, et il en est  \emph{a fortiori} de m\^eme de
l'image de $\pi_1(X)$.

\begin{rema}
La m\^eme m\'ethode permet de retrouver partiellement un r\'esultat
de Campana.
Soit $\mathcal C$ une classe de groupes qui est stable par 
sous-quotient et extension par un groupe fini. 
Supposons que deux points g\'en\'eraux d'une vari\'et\'e~$X$,
propre et normale, soient reli\'es
par une cha\^{\i}ne de sous-vari\'et\'es $Z_j$ telles que pour tout~$j$,
$\im(\pi_1(Z_j)\ra\pi_1(X))$ appartienne \`a~$\mathcal C$.
Alors, si $\mathcal C$ est de plus stable par extension,
on a $\pi_1(X)\in\mathcal C$.
Cela s'applique en particulier si $\mathcal C$ est la classe
des groupes virtuellement (pro)-polycycliques
ou (pro)-r\'esolubles.

En caract\'eristique z\'ero, Campana d\'emontre aussi dans~\cite{campana1998}
ce r\'esultat pour les classes des groupes virtuellement 
(pro)-nilpotents ou (pro)-ab\'eliens 
(qui ne sont pas stables par extension). Il serait int\'eressant de
le d\'emontrer en caract\'eristique positive.
\end{rema}

\clearpage

 \bibliographystyle{smfplain}
 \bibliography{aclab,acl}

\providecommand{\noopsort}[1]{}\providecommand{\url}[1]{\textit{#1}}
\providecommand{\bysame}{\leavevmode ---\ }
\providecommand{\og}{``}
\providecommand{\fg}{''}
\providecommand{\smfandname}{\&}
\providecommand{\smfedsname}{\'eds.}
\providecommand{\smfedname}{\'ed.}
\providecommand{\smfmastersthesisname}{M\'emoire}
\providecommand{\smfphdthesisname}{Th\`ese}
\begin{thebibliography}{10}

\bibitem{campana1992}
{\scshape F.~Campana} -- {\og Connexit\'e rationnelle des vari\'et\'es de
  {F}ano\fg}, \emph{Ann. Sci. {\'E}cole Norm. Sup.} \textbf{25} (1992), no.~5,
  p.~539--545.

\bibitem{campana1994}
\bysame , {\og Remarques sur le rev\^etement universel des vari\'et\'es
  k\"ahl\'eriennes compactes\fg}, \emph{Bull. Soc. Math. France} \textbf{122}
  (1994), no.~2, p.~255--284.

\bibitem{campana1998}
\bysame , {\og Connexit\'e ab\'elienne des vari\'et\'es k\"ahl\'eriennes
  compactes\fg}, \emph{Bull. Soc. Math. France} \textbf{126} (1998), no.~4,
  p.~483--506.

\bibitem{acl2003b}
{\scshape A.~Chambert-Loir} -- {\og Points rationnels et groupes fondamentaux:
  applications de la cohomologie $p$-adique\fg}, Séminaire Bourbaki 2002/03,
  55\textsuperscript e année, exposé 914.

\bibitem{debarre2001}
{\scshape O.~Debarre} -- \emph{Higher-dimensional algebraic geometry},
  Universitext, Springer-Verlag, New York, 2001.

\bibitem{debarre2002}
\bysame , {\og Variétés rationnellement connexes\fg}, \emph{{\upshape \`a
  paraître dans} Ast\'erisque} (2003), S\'eminaire Bourbaki, Vol.\ 2001/02,
  exposé 905.

\bibitem{deligne1972}
{\scshape P.~Deligne} -- {\og Th{\'e}orie de {H}odge {II}\fg}, \emph{Publ.
  Math. Inst. Hautes {\'E}tudes Sci.} \textbf{40} (1972), p.~5--57.

\bibitem{ekedahl1983}
{\scshape T.~Ekedahl} -- {\og Sur le groupe fondamental d'une vari\'et\'e
  unirationnelle\fg}, \emph{C. R. Acad. Sci. Paris S{\'e}r. I Math.}
  \textbf{297} (1983), no.~12, p.~627--629.

\bibitem{esnault2003}
{\scshape H.~Esnault} -- {\og Varieties over a finite field with trivial {C}how
  group of $0$-cycles have a rational point\fg}, \emph{Invent. Math.}
  \textbf{151} (2003), no.~1, p.~187--191, \url{arXiv:math.AG/0207022}.

\bibitem{grauert-r1984}
{\scshape H.~Grauert {\normalfont \smfandname} R.~Remmert} -- \emph{Coherent
  analytic sheaves}, Grundlehren der Mathematischen Wissenschaften, vol. 265,
  Springer-Verlag, Berlin, 1984.

\bibitem{dejong-s2003}
{\scshape A.~J. de~Jong {\normalfont \smfandname} J.~Starr} -- {\og Every
  rationally connected variety over the function field of a curve has a
  rational point\fg}, 2002,
  \url{http://www-math.mit.edu/~dejong/papers/familyofcurves3.dvi}.

\bibitem{kollar1995}
{\scshape J.~Koll{\'a}r} -- \emph{Shafarevich maps and automorphic forms}, M.
  B. Porter Lectures, Princeton University Press, Princeton, NJ, 1995.

\bibitem{kollar-m-m1992b}
{\scshape J.~Koll{\'a}r, Y.~Miyaoka {\normalfont \smfandname} S.~Mori} -- {\og
  Rational connectedness and boundedness of {F}ano manifolds\fg}, \emph{J.
  Differential Geom.} \textbf{36} (1992), no.~3, p.~765--779.

\bibitem{kollar-m-m1992}
\bysame , {\og Rationally connected varieties\fg}, \emph{J. Algebraic Geometry}
  \textbf{1} (1992), no.~3, p.~429--448.

\bibitem{shioda1974}
{\scshape T.~Shioda} -- {\og An example of unirational surface in
  characteristic~$p$\fg}, \emph{Math. Ann.} \textbf{211} (1974), p.~233--236.

\end{thebibliography}
\end{document}